\documentclass[11pt,letterpaper]{amsart}%
\usepackage{amsmath}
\usepackage{graphicx}%
\usepackage{amsfonts}%
\usepackage{amssymb}
\newtheorem{theorem}{Theorem}
\theoremstyle{plain}

\newtheorem{corollary}[theorem]{Corollary}

\newtheorem{lemma}[theorem]{Lemma}

\newtheorem{proposition}[theorem]{Proposition}

\begin{document}
\title[$\ell^{1}$-spreading models in mixed Tsirelson space]{$\ell^{1}$-spreading models in mixed Tsirelson space}
\author{Denny H. Leung}
\address{Department of Mathematics, National University of Singapore, 2 Science Drive
2, Singapore 117543.}
\email{matlhh@nus.edu.sg}
\author{Wee-Kee Tang}
\address{Mathematics and Mathematics Education, National Institute of Education \\
Nanyang Technological University, 1 Nanyang Walk, Singapore 637616.}
\email{wktang@nie.edu.sg}
\subjclass{}
\keywords{}

\begin{abstract}
Suppose that $(\mathcal{F}_{n})_{n=1}^{\infty}$ is a sequence of regular
families of finite subsets of $\mathbb{N}$ and $(\theta_{n})_{n=1}^{\infty}$
is a nonincreasing null sequence in $(0,1)$. The mixed Tsirelson space
$T[(\theta_{n},\mathcal{F}_{n})_{n=1}^{\infty}]$ is the completion of $c_{00}$
with respect to the implicitly defined norm
\[
\| x\| =\max\{\|x\|_{c_{0}},\sup\limits_{n}\sup\theta_{n}\sum_{i=1}^{k}\|
E_{i}x\|\} ,
\]
where the last supremum is taken over all sequences $(E_{i})_{i=1}^{k}$ in
$[\mathbb{N}]^{<\infty}$ such that $\max E_{i}<\min E_{i+1}$ and $\{\min
E_{i}:1\leq i\leq k\}\in\mathcal{F}_{n}$. Necessary and sufficient conditions
are obtained for the existence of higher order $\ell^{1}$-spreading models in
every subspace generated by a subsequence of the unit vector basis of
$T[(\theta_{n},\mathcal{F}_{n})_{n=1}^{\infty}]$.

\end{abstract}
\maketitle

\section{Preliminaries}

Mixed Tsirelson spaces were first introduced by Argyros and Deliyanni
\cite{AD}. They furnish a central class of examples in the recent development
of the structure theory of Banach spaces. In \cite{LT}, the authors computed
the Bourgain $\ell^{1}$-indices of mixed Tsirelson spaces. A stronger measure
of the finite dimensional $\ell^{1}$-structure of a Banach space is the
presence of (higher order) $\ell^{1}$-spreading models. Kutzarova and Lin
\cite{KL} showed that the Schlumprecht space \cite{S}, a fundamental example
that opened the door to much of the recent progress in the structure theory of
Banach spaces, contains an $\ell^{1}$-spreading model. Subsequently, Argyros,
Deliyanni and Manoussakis \cite{ADM} showed that if $\theta_{n+m}\geq
\theta_{n}\theta_{m}$ and $\lim_{n}\theta_{n}^{1/n}=1$, then the mixed
Tsirelson space $T[(\theta_{n},\mathcal{S}_{n})_{n=1}^{\infty}]$ contains
$\ell^{1}$-$\mathcal{S}_{\omega}$-spreading models hereditarily. In the
present paper, we consider general mixed Tsirelson spaces $T[(\theta
_{n},{\mathcal{F}}_{n})_{n=1}^{\infty}]$ and obtain necessary and sufficient
conditions for the existence of higher order $\ell^{1}$-spreading models in
every subspace generated by a subsequence of the unit vector basis.

We set the notation in the remainder of the section. Endow the power set of
$\mathbb{N}$, identified with $2^{\mathbb{N}}$, with the product topology. If
$M$ is an infinite subset of $\mathbb{N}$, denote the set of all finite,
respectively infinite, subsets of $M$ by $[M]^{<\infty},$ respectively $[M]$.
A family $\mathcal{F}\subseteq\lbrack\mathbb{N}]^{<\infty}$ is said to be
\emph{hereditary} if $G\subseteq F\in\mathcal{F}$ implies $G\in\mathcal{F}$.
It is \emph{spreading} if whenever $F=\{n_{1},\dots,n_{k}\}\in\mathcal{F}$,
$n_{1}<\dots<n_{k}$, and $m_{1}<\dots<m_{k}$ satisfy $m_{i}\geq n_{i}$, $1\leq
i\leq k$, then $\{m_{1},\dots,m_{k}\}\in\mathcal{F}$.
A \emph{regular} family is one that is hereditary, spreading and compact (as a
subset of the topological space $[\mathbb{N}]^{<\infty}$). If $E$ and $F$ are
finite subsets of $\mathbb{N}$, we write $E<F$, respectively $E\leq F$, to
mean $\max E<\min F$, respectively $\max E\leq\min F$ ($\max\emptyset=0$ and
$\min\emptyset=\infty$). We abbreviate $\{n\}<E$ and $\{n\}\leq E$ to $n<E$
and $n\leq E$ respectively. Given $\mathcal{F}\subseteq\lbrack\mathbb{N}%
]^{<\infty}$, a sequence of finite subsets $\{E_{1},\dots,E_{n}\}$ of
$\mathbb{N}$ is said to be $\mathcal{F}$-\emph{admissible} if $E_{1}%
<\dots<E_{n}$ and $\{\min E_{1},\dots,\min E_{n}\}\in\mathcal{F}$. If
$\mathcal{M}$ and $\mathcal{N}$ are regular subsets of $[\mathbb{N}]^{<\infty
}$, we let
\[
\mathcal{M}[\mathcal{N}]=\{\cup_{i=1}^{k}F_{i}:F_{i}\in\mathcal{N}\text{ for
all $i$ and }\{F_{1},\dots,F_{k}\}\text{ is }\mathcal{M}\text{-admissible}\}.
\]
Given a sequence of regular families $(\mathcal{M}_{i})$, we define
inductively $[\mathcal{M}_{1},\mathcal{M}_{2}]=\mathcal{M}_{1}[\mathcal{M}%
_{2}]$ and $[\mathcal{M}_{1},\dots,\mathcal{M}_{i+1}]=[\mathcal{M}_{1}%
,\dots,\mathcal{M}_{i}][\mathcal{M}_{i+1}]$. Also, let
\[
(\mathcal{M}_{1},\dots,\mathcal{M}_{k})=\{\cup_{i=1}^{k}M_{i}:M_{i}%
\in\mathcal{M}_{i},M_{1}<\dots<M_{k}\}.
\]
We abbreviate the $k$-fold constructions $[\mathcal{M},\dots,\mathcal{M}]$ and
$(\mathcal{M},\dots,\mathcal{M})$ as $[\mathcal{M}]^{k}$ and $(\mathcal{M}%
)^{k}$ respectively. Of primary importance are the Schreier classes as defined
in \cite{AA}. Let $\mathcal{S}_{0}=\{\{n\}:n\in\mathbb{N}\}\cup\{\emptyset\}$
and $\mathcal{S}_{1}=\{F\subseteq\mathbb{N}:|F|\leq\min F\}$. Here $|F|$
denotes the cardinality of $F$. The higher Schreier classes are defined
inductively as follows. $\mathcal{S}_{\alpha+1}=\mathcal{S}_{1}[\mathcal{S}%
_{\alpha}]$ for all $\alpha<\omega_{1}$. If $\alpha$ is a countable limit
ordinal, choose a sequence $(\alpha_{n})$ strictly increasing to $\alpha$ and
set
\[
\mathcal{S}_{\alpha}=\{F:F\in\mathcal{S}_{\alpha_{n}}\text{ for some
$n\leq\min F$}\}.
\]
It is clear that $\mathcal{S}_{\alpha}$ is a regular family for all
$\alpha<\omega_{1}$. Given a nonzero countable ordinal $\alpha$ whose Cantor
normal form is $\alpha=\omega^{\beta_{1}}\cdot m_{1}+\dots+\omega^{\beta_{n}%
}\cdot m_{n}$, we let $\mathcal{R}_{\alpha}$ be the regular family
$((\mathcal{S}_{\beta_{n}})^{m_{n}},\dots,(\mathcal{S}_{\beta_{1}})^{m_{1}})$.
If $\mathcal{F}$ is a closed subset of $[\mathbb{N}]^{<\infty}$, let
$\mathcal{F}^{\prime}$ be the set of all limit points of $\mathcal{F}$. Define
a transfinite sequence of sets $(\mathcal{F}^{(\alpha)})_{\alpha<\omega_{1}}$
as follows: $\mathcal{F}^{(0)}=\mathcal{F}$, $\mathcal{F}^{(\alpha
+1)}=(\mathcal{F}^{(\alpha)})^{\prime}$ for all $\alpha<\omega_{1}$;
$\mathcal{F}^{(\alpha)}=\cap_{\beta<\alpha}\mathcal{F}^{(\beta)}$ if $\alpha$
is a countable limit ordinal. If $\mathcal{F}$ is regular, we let
$\iota(\mathcal{F})$ be the unique ordinal $\alpha$ such that $\mathcal{F}%
^{(\alpha)}=\{\emptyset\}$. It is well known that $\iota(\mathcal{S}_{\gamma
})=\omega^{\gamma}$ for all $\gamma<\omega_{1}$ \cite[Proposition 4.10]{AA}.
Also, $\iota((\mathcal{M},\mathcal{N}))=\iota(\mathcal{N})+\iota(\mathcal{M})$
and $\iota(\mathcal{M}[\mathcal{N}])\leq\iota(\mathcal{N})\cdot\iota
(\mathcal{M})$ \cite[Proposition 10]{LT1}. In particular, $\iota
(\mathcal{R}_{\alpha})=\alpha$.

If $\mathcal{F}$ is a regular family and $K$ is a positive constant, we say
that a normalized sequence $(x_{n})$ in a Banach space is an $\ell^{1}%
$-$\mathcal{F}$-\emph{spreading model with constant} $K$ if $\|\sum_{F}%
a_{n}x_{n}\| \geq K^{-1}\sum_{F}|a_{n}|$ for all $F \in\mathcal{F}$ and all
sequences of scalars $(a_{n})$. We refer to \cite{JO} for the definitions and
in depth discussions of the $\ell^{1}$-indices $I(X)$, $I(X,K)$, $I_{b}(X)$
and $I_{b}(X,K)$ of a Banach space $X$ (assumed to have a basis in the last
two). Suffice it to say that if $X$ contains an $\ell^{1}$-$\mathcal{F}%
$-spreading model with constant $K$, then $I(X,K) \geq\iota(\mathcal{F})$.
Moreover, if the spreading model is a block basis of the basis of $X$, then
$I_{b}(X,K) \geq\iota(\mathcal{F})$.

Let $c_{00}$ be the vector space of all finitely supported real sequences and
let $(e_{k})$ be the standard unit vector basis of $c_{00}$.
For $E\in\lbrack\mathbb{N}]^{<\infty}$ and $x=\sum a_{k}e_{k}\in c_{00}$, let
$Ex=\sum_{k\in E}a_{k}e_{k}$.
Given a sequence of regular families $(\mathcal{F}_{n})_{n=1}^{\infty}$ and a
nonincreasing null sequence $(\theta_{n})_{n=1}^{\infty}$ in $(0,1)$, the
\emph{mixed Tsirelson space} $T[(\theta_{n},\mathcal{F}_{n})_{n=1}^{\infty}]$
is the completion of $c_{00}$ under the implicitly defined norm%
\begin{equation}
\Vert x\Vert=\max\{\Vert x\Vert_{c_{0}},\sup\limits_{n}\sup\theta_{n}%
\sum_{i=1}^{k}\Vert E_{i}x\Vert\}, \label{norm}%
\end{equation}
where the last supremum is taken over all $\mathcal{F}_{n}$-admissible
sequences $(E_{i})_{i=1}^{k}.$

Throughout the paper, we consider a fixed mixed Tsirelson space $X=T[(\theta
_{n},\mathcal{F}_{n})_{n=1}^{\infty}]$. Set $\alpha_{n}=\iota(\mathcal{F}%
_{n})$ for all $n$ and let $\alpha=\sup_{n}\alpha_{n}$. To avoid trivial
cases, we will assume that $\alpha_{n}>1$ for all $n$. The following
fundamental set theoretic dichotomy due to Gasparis will be used repeatedly.

\begin{theorem}
\label{G}\emph{\cite[Theorem 1.1]{G}}Let $\mathcal{F}$ and $\mathcal{G}$ be
hereditary families of finite subsets of $\mathbb{N}$ and $N$ an infinite
subset of $\mathbb{N}$. Then there exists $M\in\lbrack N]$ such that either
$\mathcal{G}\cap\left[  M\right]  ^{<\infty}\subseteq\mathcal{F}$ or
$\mathcal{F}\cap\left[  M\right]  ^{<\infty}\subseteq\mathcal{G}.$
\end{theorem}

Note that if ${\mathcal{G}}$ is a regular family, then $\iota({\mathcal{G}%
}\cap\lbrack M]^{<\infty})=\iota({\mathcal{G}})$ for all $M\in\lbrack
\mathbb{N}]$. Thus if ${\mathcal{F}}$ and ${\mathcal{G}}$ are regular families
such that $\iota({\mathcal{F}})<\iota(G)$, then for any $N\in\lbrack
\mathbb{N}]$, there exists $M\in\lbrack N]$ such that ${\mathcal{F}}%
\cap\lbrack M]^{<\infty}\subseteq{\mathcal{G}}$.

\begin{proposition}
If $\alpha=\alpha_{n}$ for some $n$ or if $\alpha$ is not of the form
$\omega^{\omega^{\xi}},$ $\xi<\omega_{1},$ then $X$
contains $\ell^{1}$-${\mathcal{R}}_{\alpha^{k}}$-spreading models hereditarily
for all $k\in\mathbb{N}.$ However, it does not contain any $\ell^{1}%
$-${\mathcal{R}}_{\alpha^{\omega}}$-spreading model.
\end{proposition}

\begin{proof}
Let $(x_{m})$ be a normalized block sequence in $X$. Under the hypothesis, for
any $k\in\mathbb{N}$, there exist $n,i\in\mathbb{N}{}$ such that $\alpha
^{k}<\alpha_{n}^{i}$. Then $\iota({\mathcal{R}}_{\alpha^{k}})<\iota
([\mathcal{F}_{n}]^{i}).$ By Theorem \ref{G} and the subsequent remark, there
exists $M\in\lbrack\mathbb{N}]^{<\infty}$ such that ${\mathcal{R}}_{\alpha
^{k}}\cap\lbrack M]^{<\infty}\subseteq\lbrack\mathcal{F}_{n}]^{i}.$ We claim
that $(x_{m})_{m\in M}$ is an $\ell^{1}$-${\mathcal{R}}_{\alpha^{k}}%
$-spreading model with constant $1/\theta^{i}$. Indeed, suppose that
$M=(m_{j})$ and $F\in{\mathcal{R}}_{\alpha^{k}},$ then $\{m_{j}:j\in
F\}\in{\mathcal{R}}_{\alpha^{k}}\cap\lbrack M]^{<\infty}\subseteq
\lbrack\mathcal{F}_{n}]^{i}.$ As a result, $\{\operatorname{supp}x_{m_{j}%
}:j\in F\}$ is $[\mathcal{F}_{n}]^{i}$-admissible$.$ Therefore, for all
$(a_{j})\in c_{00},$%
\[
\Vert\sum_{j\in F}a_{j}x_{m_{j}}\Vert\geq\theta_{n}^{i}\sum_{j\in F}\Vert
a_{j}x_{m_{j}}\Vert=\theta_{n}^{i}\sum_{j\in F}|a_{j}|.
\]

On the other hand, $I_{b}(X) =\alpha^{\omega}$ \cite[Theorem 15]{LT}. If
$\alpha\geq\omega$, then $I(X) = I_{b}(X) = \alpha^{\omega}$ by
\cite[Corollary 5.13]{JO}. By \cite[Lemma 5.11]{JO}, $I( X,K) <\alpha^{\omega
}$ for all $K\geq1.$ It follows that $X$ does not contain an $\ell^{1}%
$-${\mathcal{R}}_{\alpha^{\omega}}$-spreading model. If $\alpha< \omega$, then
$\alpha^{\omega}= \omega$ since we are assuming that $\alpha> 1$. If $(x_{n})$
is an $\ell^{1}$-$\mathcal{S}_{1}$-spreading model in $X$, then there is a
subsequence $(x_{n_{k}})$ such that $(x_{n_{2k}}-x_{n_{2k+1}})$ is equivalent
to a block basis of the unit vector basis $(e_{k})$ of $X$. It is easily
checked that $(x_{n_{2k}}-x_{n_{2k+1}})$ is an $\ell^{1}$-$\mathcal{S}_{1}%
$-spreading model. Thus $\omega\leq I_{b}(X,K)$ and hence $I_{b}(X) =
I_{b}(X,K)$, contrary to \cite[Lemma 5.7]{JO}.
\end{proof}

\section{Higher order $\ell^{1}$-spreading models}

Henceforth, we assume that $\alpha\neq\alpha_{n}$ for any $n$ and
$\alpha=\omega^{\omega^{\xi}}$ for some $0<\xi<\omega_{1}.$
For a nonzero ordinal $\alpha$ with Cantor normal form $\omega^{\beta_{1}%
}\cdot m_{1}+\dots+\omega^{\beta_{n}}\cdot m_{n}$, let $\ell(\alpha)=\beta
_{1}$. Given $m\in\mathbb{N}$ and $\varepsilon>0,$ define
\[
\gamma=\gamma(\varepsilon,m)=\max\{\ell(\alpha_{n_{s}}\dots\alpha_{n_{1}%
}):\varepsilon\theta_{n_{1}}\cdots\theta_{n_{s}}>\theta_{m}\}\text{
($\max\emptyset=0$).}%
\]
We say that the space $X$ satisfies $(\dagger)$ if%
\[%
\begin{tabular}
[c]{l}%
there exists $\varepsilon>0$ such that\ for all $\beta<\omega^{\xi},$ there
exists $m\in\mathbb{N}$\\
satisfying $\gamma(\varepsilon,m)+2+\beta<\ell(\alpha_{m})$.
\end{tabular}
\
\]
It was proved in \cite{LT} that condition $(\dagger)$ is sufficient for $X$ to
have a large $\ell^{1}$-index.

\begin{theorem}
\emph{\label{T1}\cite[Theorem 17]{LT}} Assume that $\xi\neq0.$ If $X$
satisfies \emph{(}$\dagger$\emph{)}, then $I(X) =\omega^{\omega^{\xi}\cdot2}$.
\end{theorem}

\noindent\emph{Remark.} It was shown in \cite[Corollary 18]{LT} that
($\dagger$) holds if $\xi$ is a limit ordinal.\newline 

Observe that if $X$ contains an $\ell^{1}$-$\mathcal{S}_{\omega^{\xi}}%
$-spreading model, then it actually contains $\ell^{1}$-$\mathcal{F}_{n}[
\mathcal{S}_{\omega^{\xi}}] $-spreading models for all $n.$ In this case, it
follows that $I( X) =\omega^{\omega^{\xi}\cdot2}$. Hence the next result
strengthens Theorem \ref{T1}.

\begin{theorem}
\label{T2}Suppose that $0<\xi<\omega_{1}$ and \emph{(}$\dagger$\emph{)} holds.
Then for any subsequence $( e_{n}) _{n\in M}$ of the unit vector basis $(
e_{n}) $ of $X,$ $[ ( e_{n}) _{n\in M}] $ contains an $\ell^{1}$%
-$\mathcal{S}_{\omega^{\xi}}$-spreading model.
\end{theorem}

The construction, using interlaced layers of vectors of differing
complexities, is based on the method pioneered by Kutzarova and Lin
(\cite{KL}) and subsequently refined and extended by Argyros
et.\ al.\ (\cite{ADKM}). As in \cite{LT}, we calculate the norms of vectors in
$X$ by means of admissible trees. Let us recall the relevant procedure and set
the notation. A \emph{tree} in $[ \mathbb{N}] ^{<\infty}$ is a finite
collection of elements $( E_{i}^{m}) ,$ $0\leq m\leq r,$ $1\leq i\leq k( m) ,$
in $[ \mathbb{N}] ^{<\infty}$ so that for each $m,$ $E_{1}^{m}<E_{2}^{m}%
<\dots<E_{k( m) }^{m},$ and that every $E_{i}^{m+1}$ is a subset of some
$E_{j}^{m}.$ The elements $E_{i}^{m}$ are called \emph{nodes} of the tree. Any
node $E_{i}^{m}$ is said to be of \emph{level} $m.$ Nodes at level $0$ are
called \textit{roots. }If $E_{i}^{n}\subseteq E_{j}^{m}$ and $n>m,$ we say
that $E_{i}^{n}$ is a \emph{descendant }of $E_{j}^{m}$ and $E_{j}^{m}$ is an
\emph{ancestor} of $E_{i}^{n}.$ If, in the above notation, $n=m+1,$ then
$E_{i}^{n}$ is said to be an \emph{immediate successor} of $E_{j}^{m},$ and
$E_{j}^{m}$ the \emph{immediate predecessor }of $E_{i}^{n}.$ Nodes with no
descendants are called \emph{terminal nodes} or \emph{leaves} of the tree. The
set of all leaves of a tree $\mathcal{T}$ is denoted by $\mathcal{L}(
\mathcal{T}) .$ A tree $( E_{i}^{m}) ,$ $0\leq m<r,$ $1\leq i\leq k( m) ,$ is
$( \mathcal{F}_{n}) $-admissible if $k( 0) =1$ and for every $m$ and $i,$ the
collection $( E_{j}^{m+1}) $ of all immediate successors of $E_{i}^{m}$ is an
$\mathcal{F}_{n}$-admissible collection for some $n\in\mathbb{N}.$ Given an $(
\mathcal{F}_{n}) $-admissible tree $( E_{i}^{m}) ,$ we define the
\emph{history} of the individual nodes inductively as follows. Let $h(
E_{1}^{0}) =( 0) .$ If $h( E_{i}^{m}) $ has been defined and the collection $(
E_{j}^{m+1}) $ of all immediate successors of $E_{i}^{m}$ forms an
$\mathcal{F}_{n}$-admissible collection, then define $h( E_{j}^{m+1}) $ to be
the $(m+2)$-tuple $( h( E_{i}^{m}) ,n)$.
Finally, assign $( ( \theta_{n}) \text{-compatible}) $ \emph{tags} to the
nodes by defining $t( E_{i}^{m}) =\prod_{j=0}^{m}\theta_{n_{j}}$ if $h(
E_{i}^{m}) =( n_{0},n_{1},\dots,n_{m}) $ $( \theta_{0}=1) .$ If $x\in c_{00}$
and $\mathcal{T}$ is an $( \mathcal{F}_{n}) $-admissible tree, let
$\mathcal{T}x=\sum t( E) \| Ex\| _{c_{0}},$ where the sum is taken over all
leaves in $\mathcal{T}.$ It is easily observed that $\| x\| =\max\{
\mathcal{T}x:\mathcal{T}\text{ is an }( \mathcal{F}_{n}) \text{-admissible
tree}\} .$

We are now ready to set up for the main step of the calculation. Let
$\varepsilon\in(0,1)$ be given. For $r\in\mathbb{N},$ let $\mathcal{N}%
_{r}=\{(0,n_{1},...,n_{s}):\varepsilon\theta_{n_{1}}\cdots\theta_{n_{s}%
}>\theta_{r}\}.$ Then $\gamma(\varepsilon,r)=\max\{\ell(\alpha_{n_{s}%
}...\alpha_{n_{1}}):(0,n_{1},...,n_{s})\in\mathcal{N}_{r}\}.$ Assume
$\delta\in(0,1),$ $p<q$ and $\eta$ are given such that $\gamma(\varepsilon
,p)<\eta<\omega^{\xi}.$ Let
\[
K_{\delta,p,\eta}=\{(0,n_{1},...,n_{s}):\theta_{n_{1}}\cdots\theta_{n_{s}%
}>\delta\theta_{p},\text{ }\ell(\alpha_{n_{s}}...\alpha_{n_{1}})<\eta\}.
\]
Also assume that $M\in\lbrack\mathbb{N}]$ satisfies $[\mathcal{F}_{n_{1}%
},...,\mathcal{F}_{n_{s}}]\cap\lbrack M]^{<\infty}\subseteq\mathcal{S}_{\eta}$
whenever $(0,n_{1},...,n_{s})\in K_{\delta,p,\eta}.$ Suppose that vectors
$x_{1}$ and $x_{2}$ are given so that $x_{1}={\theta_{p}^{-1}}%
{\displaystyle\sum\limits_{i=1}^{r}}
{a_{i}}e_{m_{i}}$, $x_{2}=%
{\displaystyle\sum\limits_{i=1}^{r}}
a_{i}z_{i},$ $x=x_{1}+x_{2},$ and
\begin{align*}
\Vert x_{1}\Vert_{\mathcal{S}_{\eta}}  &  \leq\dfrac{\delta}{|K_{\delta
,p,\eta}|+1},\\
\{m_{1},m_{2},...,m_{r}\}  &  \in\mathcal{S}_{\eta+1}\cap\lbrack M]^{<\infty
},\\
\Vert x_{1}\Vert_{\ell^{1}}  &  =\dfrac{1}{\theta_{p}},\\
m_{1}  &  <z_{1}<...<m_{r}<z_{r}.
\end{align*}
If $y=\sum a_{k}e_{k}\in c_{00}$ and ${\mathcal{F}}$ is a regular family, let
$\Vert y\Vert_{{\mathcal{F}}}=\sup_{F\in{\mathcal{F}}}\sum_{k\in F}|a_{k}|$.

\begin{proposition}
\label{P4} Let $x$ be given as above. For any admissible tree $\mathcal{T}$,
there exist an admissible tree $\mathcal{T}^{\prime}$ and disjoint sets
$J_{1}$ and $J_{2}$ such that

(1) $\mathcal{T}^{\prime}$ is $( p,q) $-restricted, i.e., for all
$E\in\mathcal{L}( \mathcal{T}^{\prime}) ,$ there exists $G\in\mathcal{T}%
^{\prime}$ containing $E$ such that $h( G) \in\mathcal{N}_{q}\setminus
\mathcal{N}_{p},$

(2) $\mathcal{T}x\leq\mathcal{T}(
{\displaystyle\sum\limits_{i\in J_{1}}}
\dfrac{a_{i}}{\theta_{p}}e_{m_{i}}+%
{\displaystyle\sum\limits_{i\in J_{2}}}
a_{i}z_{i}) +\mathcal{T}^{\prime}x_{2}+\delta+\dfrac{\theta_{q}}%
{\varepsilon\theta_{p}}.$
\end{proposition}

\begin{proof}
Choose $m_{r+1}>\max\operatorname*{supp}z_{r}.$ We may assume without loss of
generality that the root of $\mathcal{T}$ is the integer interval
$[m_{1},m_{r+1}]$, that every node in $\mathcal{T}$ is an integer interval,
and that every leaf in $\mathcal{T}$ is a singleton. For each $i\leq r,$ let
$\mathcal{E}_{i}=\{E\in\mathcal{L}(\mathcal{T}):E\subseteq\operatorname*{supp}%
z_{i}\}.$ Define%
\begin{align*}
I_{1}  &  =\{i:\mathcal{E}_{i}\neq\emptyset,\text{ }\{m_{i}\}\in
\mathcal{L}(\mathcal{T})\},\\
I_{2}  &  =\{i:\mathcal{E}_{i}\neq\emptyset,\text{ }\{m_{i}\}\notin
\mathcal{L}(\mathcal{T})\},\text{ and}\\
I_{3}  &  =\{i:\mathcal{E}_{i}=\emptyset,\text{ }\{m_{i}\}\in\mathcal{L}%
(\mathcal{T})\}.
\end{align*}
If $\{m_{i}\}\in\mathcal{L}(\mathcal{T}),$ we write $t_{i}$ for the tag
$t(\{m_{i}\}).$ Observe that
\begin{align}
\mathcal{T}x  &  =\sum_{E\in\mathcal{L}(\mathcal{T})}t(E)\Vert Ex\Vert_{c_{0}%
}\label{D1}\\
&  \leq\sum_{i\in I_{1}\cup I_{3}}t_{i}\dfrac{|a_{i}|}{\theta_{p}}+\sum_{i\in
I_{1}\cup I_{2}}\sum_{E\in\mathcal{E}_{i}}|a_{i}|t(E)\Vert Ez_{i}\Vert_{c_{0}%
}.\nonumber
\end{align}
For each $i\in I_{1},$ let $F_{i}$ be the smallest (by set inclusion) node in
$\mathcal{T}$ such that $\{m_{i},m_{i+1}\}\subseteq F_{i}$, then let $G_{i}$
be the immediate successor of $F_{i}$ containing $m_{i}.$ Note that if
$i_{1},\,i_{2}\in I_{1}$ and $i_{1}<i_{2},$ then $G_{i_{1}}\neq G_{i_{2}}.$
For otherwise, since $G_{i_{1}}=G_{i_{2}}$ is an integer interval,
$\{m_{i_{1}},m_{i_{1}+1}\}\subseteq G_{i_{1}}\subsetneqq F_{i_{1}},$ contrary
to the choice of $F_{i_{1}}.$ Subdivide $I_{1}$ into $I_{1}^{\prime}%
,I_{1}^{\prime\prime},$ and $I_{1}^{\prime\prime\prime}$ according to whether
$h(G_{i})\in\mathcal{N}_{p},$ $h(G_{i})\in\mathcal{N}_{q}\setminus
\mathcal{N}_{p},$ or $h(G_{i})\notin\mathcal{N}_{q}.$ Suppose $i\in
I_{1}^{\prime}.$ Then $h(G_{i})=(0,n_{1},...,n_{s})\in\mathcal{N}_{p}.$ It
follows that $\theta_{n_{1}}\cdots\theta_{n_{s}}>\delta\theta_{p}$ and
$\ell(\alpha_{n_{s}}...\alpha_{n_{1}})\leq\gamma(\varepsilon,p)<\eta.$ Thus
$h(G_{i})\in K_{\delta,p,\eta}.$ Hence
\begin{align}
\sum_{i\in I_{1}^{\prime}}t_{i}\dfrac{|a_{i}|}{\theta_{p}}  &  \leq\sum_{i\in
I_{1}^{\prime}}t(G_{i})\Vert G_{i}x_{1}\Vert_{c_{0}}\label{D2}\\
&  \leq\sum_{(0,n_{1},...,n_{s})\in K_{\delta,p,\eta}}\sum_{h(G)=(0,n_{1}%
,...,n_{s})}t(G)\Vert Gx_{1}\Vert_{c_{0}}\nonumber\\
&  \leq\sum_{(0,n_{1},...,n_{s})\in K_{\delta,p,\eta}}\Vert x_{1}%
\Vert_{\lbrack\mathcal{F}_{n_{1}},...,\mathcal{F}_{n_{s}}]}\nonumber\\
&  \leq|K_{\delta,p,\eta}|\Vert x_{1}\Vert_{\mathcal{S}_{\eta}}<\delta
.\nonumber
\end{align}
The next to last inequality holds since for any $(0,n_{1},\dots,n_{s})$, the
set $\{G\in{\mathcal{T}}:h(G)=(0,n_{1},\dots,n_{s})\}$ is $[{\mathcal{F}%
}_{n_{1}},\dots,{\mathcal{F}}_{n_{s}}]$-admissible. Also,%

\begin{equation}
\sum_{i\in I_{1}^{\prime\prime\prime}}t_{i}\dfrac{| a_{i}| }{\theta_{p}}%
\leq\sum_{i\in I_{1}^{\prime\prime\prime}}t( G_{i}) \dfrac{| a_{i}| }%
{\theta_{p}}\leq\dfrac{\theta_{q}}{\varepsilon}\| x_{1}\| _{\ell^{1}}%
=\dfrac{\theta_{q}}{\varepsilon\theta_{p}}. \label{D3}%
\end{equation}
Define $J_{1}=I_{1}^{\prime\prime}\cup I_{3}$, $J_{2}=I_{1}^{\prime}\cup
I_{1}^{\prime\prime\prime}\cup I_{2}$ and let $\mathcal{T}^{\prime}$ be the
subtree of $\mathcal{T}$ consisting of all nodes in $\bigcup_{i\in
I_{1}^{\prime\prime}}\mathcal{E}_{i}$ together with their ancestors. Clearly
$J_{1}$ is disjoint from $J_{2}.$ Note that if $E\in\mathcal{L}(
\mathcal{T}^{\prime}) ,$ then $E\in\mathcal{E}_{i}$ for some $i\in
I_{1}^{\prime\prime}.$ Since $m_{i}<E<m_{i+1}$ and $m_{i},m_{i+1}$ are both
contained in the integer interval $F_{i},$ $E\varsubsetneqq F_{i}.$ Hence
there exists an immediate successor $H$ of $F_{i}$ such that $E\subseteq H.$
But $h( H) =h( G_{i}) $ as $H$ and $G_{i}$ are both immediate successors of
$F_{i}.$ Thus $h( H) \in\mathcal{N}_{q}\setminus\mathcal{N}_{p}.$ This shows
that $\mathcal{T}^{\prime}$ is $( p,q) $-restricted. Applying $(
\text{\ref{D2}}) $ and $( \text{\ref{D3}}) $ to (\ref{D1}), we see that%
\begin{align*}
\mathcal{T}x  &  \leq\delta+\dfrac{\theta_{q}}{\varepsilon\theta_{p}}%
+\sum_{i\in I_{1}^{\prime\prime}\cup I_{3}}t_{i}\dfrac{| a_{i}| }{\theta_{p}%
}+\sum_{i\in I_{1}\cup I_{2}}\sum_{E\in\mathcal{E}_{i}}| a_{i}| t( E) \|
Ez_{i}\| _{c_{0}}\\
&  =\delta+\dfrac{\theta_{q}}{\varepsilon\theta_{p}}+\sum_{i\in J_{1}}%
t_{i}\dfrac{| a_{i}| }{\theta_{p}}+( \sum_{i\in J_{2}}+\sum_{i\in
I_{1}^{\prime\prime}}) ( \sum_{E\in\mathcal{E}_{i}}| a_{i}| t( E) \| Ez_{i}\|
_{c_{0}})\\
&  \leq\delta+\dfrac{\theta_{q}}{\varepsilon\theta_{p}}+\mathcal{T}(
\sum_{i\in J_{1}}\dfrac{a_{i}}{\theta_{p}}e_{m_{i}}+\sum_{i\in J_{2}}%
a_{i}z_{i}) +\mathcal{T}^{\prime}( \sum_{i=1}^{r}a_{i}z_{i}) ,
\end{align*}
as required.
\end{proof}

Assume that $X$ satisfies $(\dagger).$ The next step is to iterate the
construction in Proposition \ref{P4} to generate vectors with an arbitrary
number of layers. The key observation is that these vectors are uniformly
bounded. The corresponding layers in the vectors will interact to give the
desired finite dimensional $\ell^{1}$ behavior. Let $\varepsilon$ be the
constant given by condition $(\dagger).$ Suppose $(\beta_{n})$ is the sequence
of ordinals increasing to $\omega^{\xi}$ that defines $\mathcal{S}%
_{\omega^{\xi}}.$ Given any $M_{0}\in\lbrack\mathbb{N}],$ we choose sequences
$(p_{n})$, $(q_{n})$ in $\mathbb{N},$ a decreasing sequence of infinite
subsets $(M_{n})$ of $M_{0}$ and a sequence of countable ordinals $(\eta_{n})$
less than $\omega^{\xi}$ in the following manner. Pick $p_{1}\in\mathbb{N}$ so
that $\theta_{p_{1}}\leq{\varepsilon^{2}}/{4}$ and $\gamma(\varepsilon
,p_{1})+2+\beta_{1}<\ell(\alpha_{p_{1}}).$ Define $\eta_{1}=$ $\gamma
(\varepsilon,p_{1})+1.$ Then choose $q_{1}\in\mathbb{N}$ so that
$\theta_{q_{1}}\leq{\varepsilon\theta_{p_{1}}}/{4}.$ Since $\eta_{1}%
+1+\beta_{1}<\ell(\alpha_{p_{1}})$ and $\ell(\alpha_{n_{s}}\cdots\alpha
_{n_{1}})<\eta_{1}$ for all $(0,n_{1},...,n_{s})\in K_{4^{-1},p_{1}\eta_{1}},$
by the remark following Theorem \ref{G}, there exists $M_{1}\in\lbrack M_{0}]$
such that $\mathcal{S}_{\beta_{1}}[\mathcal{S}_{\eta_{1}+1}]\cap\lbrack
M_{1}]^{<\infty}\subseteq\mathcal{F}_{p_{1}}$ and $[\mathcal{F}_{n_{1}%
},...,\mathcal{F}_{n_{s}}]\cap\lbrack M_{1}]^{<\infty}\subseteq\mathcal{S}%
_{\eta_{1}}$whenever $(0,n_{1},...,n_{s})\in K_{4^{-1},p_{1},\eta_{1}}.$
Assume that the sequences have been chosen up to $n-1.$ Pick $p_{n}>q_{n-1}$
so that $\theta_{p_{n}}\leq{\varepsilon^{2}}/{4^{n}}$ and
\[
\gamma(\varepsilon,p_{n})+2+\gamma(\varepsilon,q_{n-1})+2+\eta_{n-1}%
+1+...+\eta_{1}+1+\beta_{n}<\ell(\alpha_{p_{n}}).
\]
Define $\eta_{n}=\gamma(\varepsilon,p_{n})+\gamma(\varepsilon,q_{n-1})+1.$
Then choose $q_{n}>p_{n}$ so that $\theta_{q_{n}}\leq{\varepsilon\theta
_{p_{n}}}/{4^{n}}.$ Since $\eta_{n}+1+...+\eta_{1}+1+\beta_{n}<\ell
(\alpha_{p_{n}})$ and $\ell(\alpha_{n_{s}}...\alpha_{n_{1}})<\eta_{n}$ if
$(0,n_{1},...n_{s})\in K_{4^{-n},p_{n},\eta_{n}},$ there exists $M_{n}%
\in\lbrack M_{n-1}]$ so that
\[
\mathcal{S}_{\beta_{n}}[\mathcal{S}_{\eta_{1}+1},...,\mathcal{S}_{\eta_{n}%
+1}]\cap\lbrack M_{n}]^{<\infty}\subseteq\mathcal{F}_{p_{n}}%
\]
and $[\mathcal{F}_{n_{1}},...,\mathcal{F}_{n_{s}}]\cap\lbrack M_{n}]^{<\infty
}\subseteq\mathcal{S}_{\eta_{n}}$ if $(0,n_{1},...,n_{s})\in K_{4^{-n}%
,p_{n},\eta_{n}}.$ This completes the inductive construction. For every $n$,
let $Z(p_{n})$ be the set of all vectors $x$ in $c_{00}$ such that $\Vert
x\Vert_{\ell^{1}}=\theta_{p_{n}}^{-1},$ $\operatorname*{supp}x\in
\mathcal{S}_{\eta_{n}+1}\cap\lbrack M_{n}]^{<\infty}$ and $\Vert
x\Vert_{\mathcal{S}_{\eta_{n}}}\leq4^{-n}(|K_{4^{-n},p_{n}\eta_{n}}|+1)^{-1}.$
The set $Z(p_{n})$ is nonempty by Proposition 3.6 in \cite{OTW}. Inductively,
for $n,k\in\mathbb{N},$ let $Z(p_{n},p_{n+1},...,p_{n+k})$ consists of all
vectors of the form $\theta_{p_{n}}^{-1}\sum_{i=1}^{r}a_{i}e_{m_{i}}%
+\sum_{i=1}^{r}a_{i}z_{i},$ where $m_{1}<z_{1}<...<m_{r}<z_{r},$
$\theta_{p_{n}}^{-1}\sum_{i=1}^{r}a_{i}e_{m_{i}}\in Z(p_{n})$ and $z_{i}\in
Z(p_{n+1},...,p_{n+k}),$ $1\leq i\leq r.$ Recall that an admissible tree
$\mathcal{T}$ is said to be $(p,q)$-restricted if every leaf $E\in\mathcal{T}$
is contained in some node $G\in\mathcal{T}$ with $h(G)\in\mathcal{N}%
_{q}\setminus\mathcal{N}_{p}.$ In the following, a $(p_{0},q_{0})$-restricted
tree is one without any restriction placed on it.

\begin{lemma}
\label{L5}Let $x$ be a vector finitely supported in $M_{n}$ and suppose that
$\| x\| _{\mathcal{S}_{\eta_{n}}}\leq4^{-n}( | K_{4^{-n},p_{n}\eta_{n}}| +1)
^{-1}$. If $0\leq m<n$ and $\mathcal{T}$ is a $( p_{m},q_{m}) $-restricted
admissible tree, then%
\[
\mathcal{T}x\leq\biggl\{
\begin{tabular}
[c]{ll}%
$4^{-n}+\dfrac{\theta_{p_{n}}}{\varepsilon}\| x\| _{\ell^{1}}$ & if $m=0$\\
$4^{-n}+\theta_{p_{n}}\| x\| _{\ell^{1}}( 4^{-n}+4^{-m}) $ & if $0<m<n.$%
\end{tabular}
\]
\end{lemma}

\begin{proof}
First assume that $m=0.$ Observe that $\mathcal{N}_{p_{n}}\subseteq
K_{4^{-n},p_{n}\eta_{n}}.$ Indeed, if $(0,n_{1},...,n_{s})\in\mathcal{N}%
_{p_{n}},$ then $\ell(\alpha_{n_{s}}...\alpha_{n_{1}})\leq\gamma
(\varepsilon,p_{n})<\eta_{n}$ and $\theta_{n_{1}}\cdots\theta_{n_{s}}%
>\theta_{p_{n}}/\varepsilon>4^{-n}\theta_{p_{n}}.$ Thus $(0,n_{1}%
,...,n_{s})\in K_{4^{-n},p_{n}\eta_{n}}.$ For a fixed $(0,n_{1},...,n_{s}),$
$\{E\in\mathcal{L}(\mathcal{T}):h(E)\in(0,n_{1},...,n_{s})\}$ is
$[\mathcal{F}_{n_{1}},...,\mathcal{F}_{n_{s}}]$-admissible. Hence if
$(0,n_{1},...,n_{s})\in\mathcal{N}_{p_{n}}\subseteq K_{4^{-n},p_{n}\eta_{n}},$
then%
\[
\sum_{\substack{E\in\mathcal{L}(\mathcal{T})\\h(E)=(0,n_{1},...,n_{s}%
)}}t(E)\Vert Ex\Vert_{c_{0}}\leq\Vert x\Vert_{\lbrack\mathcal{F}_{n_{1}%
},...,\mathcal{F}_{n_{s}}]\cap\lbrack M_{n}]^{<\infty}}\leq\Vert
x\Vert_{\mathcal{S}_{\eta_{n}}}.
\]
Therefore,%
\begin{align*}
\mathcal{T}x  &  \leq\biggl(\sum_{\substack{E\in\mathcal{L}(\mathcal{T}%
)\\h(E)\in\mathcal{N}_{p_{n}}}}+%
{\displaystyle\sum_{\substack{E\in\mathcal{L}(\mathcal{T})\\h(E)\notin
\mathcal{N}_{p_{n}}}}}
\biggr)t(E)\Vert Ex\Vert_{c_{0}}\\
&  \leq\sum_{(0,n_{1},...,n_{s})\in K_{4^{-n},p_{n}\eta_{n}}}\sum
_{\substack{E\in\mathcal{L}(\mathcal{T})\\h(E)=(0,n_{1},...,n_{s})}}t(E)\Vert
Ex\Vert_{c_{0}}\\
&  +\sum_{\substack{E\in\mathcal{L}(\mathcal{T})\\h(E)\notin\mathcal{N}%
_{p_{n}}}}\dfrac{\theta_{p_{n}}}{\varepsilon}\Vert Ex\Vert_{c_{0}}\\
&  \leq|K_{4^{-n},p_{n}\eta_{n}}|\Vert x\Vert_{\mathcal{S}_{\eta_{n}}}%
+\dfrac{\theta_{p_{n}}}{\varepsilon}\Vert x\Vert_{\ell^{1}}\\
&  \leq4^{-n}+\dfrac{\theta_{p_{n}}}{\varepsilon}\Vert x\Vert_{\ell^{1}}.
\end{align*}

Assume that $0<m<n$. If $E\in\mathcal{L}(\mathcal{T}),$ pick $G\in\mathcal{T}$
so that $E\subseteq G$ and $h(G)\in\mathcal{N}_{q_{m}}\setminus\mathcal{N}%
_{p_{m}}.$ Write $h(G)=(0,n_{1},...,n_{s})$ and $h(E)=(0,n_{1},...,n_{t}),$
$t\geq s.$ If $(0,n_{s+1},...,n_{t})\in\mathcal{N}_{p_{n}},$ then $\ell
(\alpha_{n_{t}}...\alpha_{n_{s+1}})\leq\gamma(\varepsilon,p_{n}).$ Since
$h(G)\in\mathcal{N}_{q_{m}}\subseteq\mathcal{N}_{q_{n-1}}$, we also have
$\ell(\alpha_{n_{s}}...\alpha_{n_{1}})\leq\gamma(\varepsilon,q_{n-1})$.
Therefore,%
\begin{align*}
\ell(\alpha_{n_{t}}...\alpha_{n_{s+1}}\alpha_{n_{s}}...\alpha_{n_{1}})  &
=\ell(\alpha_{n_{t}}...\alpha_{n_{s+1}})+\ell(\alpha_{n_{s}}...\alpha_{n_{1}%
})\\
&  \leq\gamma(\varepsilon,p_{n})+\gamma(\varepsilon,q_{n-1})<\eta_{n}.
\end{align*}
It follows that if $(0,n_{s+1},...,n_{t})\in\mathcal{N}_{p_{n}}$ and
$t(E)>4^{-n}\theta_{p_{n}},$ then $h(E)\in K_{4^{-n},p_{n}\eta_{n}}.$ Thus,%
\begin{align}
&  \sum_{\substack{E\in\mathcal{L}(\mathcal{T})\\(0,n_{s+1},...,n_{t}%
)\in\mathcal{N}_{p_{n}}}}t(E)\Vert Ex\Vert_{c_{0}}\label{D4}\\
&  \leq\sum_{\substack{E\in\mathcal{L}(\mathcal{T})\\h(E)\in K_{4^{-n}%
,p_{n},\eta_{n}}}}t(E)\Vert Ex\Vert_{c_{0}}+\sum_{\substack{E\in
\mathcal{L}(\mathcal{T})\\t(E)\leq4^{-n}\theta_{p_{n}}}}t(E)\Vert
Ex\Vert_{c_{0}}\nonumber\\
&  \leq|K_{4^{-n},p_{n},\eta_{n}}|\Vert x\Vert_{\mathcal{S}_{\eta_{n}}}%
+4^{-n}\theta_{p_{n}}\Vert x\Vert_{\ell^{1}}\nonumber\\
&  \leq4^{-n}+4^{-n}\theta_{p_{n}}\Vert x\Vert_{\ell^{1}}.\nonumber
\end{align}

On the other hand, if $(0,n_{s+1},...,n_{t})\notin\mathcal{N}_{p_{n}},$ then
$\varepsilon\theta_{n_{s+1}}\cdots\theta_{n_{t}}\leq\theta_{p_{n}}.$
Similarly, $\varepsilon\theta_{n_{1}}\cdots\theta_{n_{s}}\leq\theta_{p_{m}}$
since $h(G)\notin\mathcal{N}_{p_{m}}.$ Hence $t(E)=\theta_{n_{1}}\cdots
\theta_{n_{s}}\theta_{n_{s+1}}\cdots\theta_{n_{t}}\leq\theta_{p_{m}}%
\theta_{p_{n}}/\varepsilon^{2}.$ Thus
\begin{equation}
\sum_{\substack{E\in\mathcal{L}(\mathcal{T})\\(0,n_{s+1},...,n_{t}%
)\notin\mathcal{N}_{p_{n}}}}t(E)\Vert Ex\Vert_{c_{0}}\leq\dfrac{\theta_{p_{m}%
}\theta_{p_{n}}}{\varepsilon^{2}}\Vert x\Vert_{\ell^{1}}\leq4^{-m}%
\theta_{p_{n}}\Vert x\Vert_{\ell^{1}}. \label{D5}%
\end{equation}
Combining (\ref{D4}) and (\ref{D5}) completes the proof.
\end{proof}

\begin{lemma}
\label{L6}Let $x$ be a vector in $Z(p_{n},...,p_{n+k}) ,$ where $n\in
\mathbb{N}$ and $k\in\mathbb{N\cup}\{ 0\} .$ If $0\leq m<n$ and $\mathcal{T}$
is a $( p_{m},q_{m}) $-restricted admissible tree$,$ then%
\[
\mathcal{T}x\leq4^{-( n-1) }\sum_{j=0}^{k}2^{-j}+\biggl\{
\begin{tabular}
[c]{ll}%
$\dfrac{1}{\varepsilon}-3\cdot4^{-( n+k) }$ & if $m=0$\\
$4^{-m}$ & if $0<m<n.$%
\end{tabular}
\]
\end{lemma}

\begin{proof}
Observe that any vector $x\in Z( p_{n}) $ satisfies the hypothesis of Lemma
\ref{L5} and that $\| x\| _{\ell^{1}}=\theta_{p_{n}}^{-1}.$ The result for
$k=0$ follows from the same lemma.

Now suppose the result holds for some $k$ and consider a vector $x\in Z(
p_{n},...,p_{n+k+1})$ and a $(p_{m},q_{m})$-restricted admissible tree
$\mathcal{T}$, $0 \leq m < n$. Write $x=\theta_{p_{n}}^{-1}\sum_{i=1}^{r}%
a_{i}e_{m_{i}}+\sum_{i=1}^{r}a_{i}z_{i}=x_{1}+x_{2}$ according to the
definition of $Z(p_{n},...,p_{n+k+1}).$ One can easily verify all the
conditions preceding Proposition \ref{P4} with the parameters $\delta=4^{-n},$
$p=p_{n},$ $q=q_{n},$ $M=M_{n},$ and $\eta=\eta_{n}.$ By Proposition \ref{P4},
we obtain a $( p_{n},q_{n}) $-restricted admissible tree $\mathcal{T}^{\prime
}$ and disjoint sets $J_{1}$ and $J_{2}$ so that%
\begin{align*}
\mathcal{T}x  &  \leq\mathcal{T}\bigl(
{\displaystyle\sum\limits_{i\in J_{1}}}
\dfrac{a_{i}}{\theta_{p_{n}}}e_{m_{i}}+%
{\displaystyle\sum\limits_{i\in J_{2}}}
a_{i}z_{i}\bigr) +\mathcal{T}^{\prime}x_{2}+4^{-n}+\dfrac{\theta_{q_{n}}%
}{\varepsilon\theta_{p_{n}}}\\
&  \leq\mathcal{T}\bigl(
{\displaystyle\sum\limits_{i\in J_{1}}}
\dfrac{a_{i}}{\theta_{p_{n}}}e_{m_{i}}+%
{\displaystyle\sum\limits_{i\in J_{2}}}
a_{i}z_{i}\bigr) +\mathcal{T}^{\prime}x_{2}+2\cdot4^{-n}.
\end{align*}
By Lemma \ref{L5},%
\[
\mathcal{T}\bigl(
{\displaystyle\sum\limits_{i\in J_{1}}}
\dfrac{a_{i}}{\theta_{p_{n}}}e_{m_{i}}\bigr)  \leq\biggl\{
\begin{tabular}
[c]{ll}%
$4^{-n}+\dfrac{1}{\varepsilon}\sum_{i\in J_{1}}| a_{i}| $ & if $m=0$\\
$4^{-n}+(4^{-n}+4^{-m}) \sum_{i\in J_{1}}|a_{i}| $ & if $m\neq0.$%
\end{tabular}
\]
Moreover, by the inductive hypothesis,%

\begin{align*}
\mathcal{T}\bigl(
{\displaystyle\sum\limits_{i\in J_{2}}}
a_{i}z_{i}\bigr)  &  \leq%
{\displaystyle\sum\limits_{i\in J_{2}}}
| a_{i}| \sup_{i\in J_{2}}\mathcal{T}z_{i}\\
&  \leq%
{\displaystyle\sum\limits_{i\in J_{2}}}
| a_{i}| \biggl(4^{-n}\sum_{j=0}^{k}2^{-j}+\biggl\{
\begin{tabular}
[c]{ll}%
$\dfrac{1}{\varepsilon}-3\cdot4^{-( n+k+1) }$ & if $m=0$\\
$4^{-m}$ & if $m\neq0$.
\end{tabular}
\biggr)%
\end{align*}
Using the fact that
\[
u\sum_{i\in J_{1}}| a_{i}| +v%
{\displaystyle\sum\limits_{i\in J_{2}}}
| a_{i}| \leq\max\{ u,v\} \sum_{i\in J_{1}\cup J_{2}}| a_{i}| \leq\max\{
u,v\}
\]
if $u,v\geq0,$ we see that
\begin{align*}
\mathcal{T}\biggl( {\displaystyle\sum\limits_{i\in J_{1}}} \dfrac{a_{i}%
}{\theta_{p_{n}}}  &  e_{m_{i}}+ {\displaystyle\sum\limits_{i\in J_{2}}}
a_{i}z_{i}\biggr)\\
\leq4^{-n}  &  +\sum_{i\in J_{1}}|a_{i}|\biggl(\biggl\{
\begin{tabular}
[c]{ll}%
$\dfrac{1}{\varepsilon}$ & if $m=0$\\
$4^{-n}+4^{-m}$ & if $m\neq0$%
\end{tabular}
\biggr)\\
&  +%
{\displaystyle\sum\limits_{i\in J_{2}}}
| a_{i}| \biggl( 4^{-n}\sum_{j=0}^{k}2^{-j}+\biggl\{
\begin{tabular}
[c]{ll}%
$\dfrac{1}{\varepsilon}-3\cdot4^{-( n+k+1) }$ & if $m=0$\\
$4^{-m}$ & if $m\neq0$%
\end{tabular}
\biggr)\\
\leq4^{-n}  &  +4^{-n}\sum_{j=0}^{k}2^{-j}+\biggl\{
\begin{tabular}
[c]{ll}%
$\dfrac{1}{\varepsilon}-3\cdot4^{-( n+k+1) }$ & if $m=0$\\
$4^{-m}$ & if $m\neq0.$%
\end{tabular}
\end{align*}
Since $\mathcal{T}^{\prime}$ is $( p_{n},q_{n}) $-restricted, the inductive
hypothesis yields
\[
\mathcal{T}^{\prime}x_{2}\leq4^{-n}\sum_{j=0}^{k}2^{-j}+4^{-n}.
\]
Therefore,%
\begin{align*}
\mathcal{T}x  &  \leq4^{-n}+4^{-n}\sum_{j=0}^{k}2^{-j}+\biggl\{
\begin{tabular}
[c]{ll}%
$\dfrac{1}{\varepsilon}-3\cdot4^{-( n+k+1) }$ & if $m=0$\\
$4^{-m}$ & if $m\neq0$%
\end{tabular}
\\
&  \quad\quad\quad+4^{-n}\sum_{j=0}^{k}2^{-j}+4^{-n}+2\cdot4^{-n}\\
&  =4\cdot4^{-n}+2\cdot4^{-n}\sum_{j=0}^{k}2^{-j}+\biggl\{
\begin{tabular}
[c]{ll}%
$\dfrac{1}{\varepsilon}-3\cdot4^{-( n+k+1) }$ & if $m=0$\\
$4^{-m}$ & if $m\neq0$%
\end{tabular}
\\
&  =4^{-( n-1) }+4^{-( n-1) }\sum_{j=0}^{k}2^{-( j+1) }+\biggl\{
\begin{tabular}
[c]{ll}%
$\dfrac{1}{\varepsilon}-3\cdot4^{-( n+k+1) }$ & if $m=0$\\
$4^{-m}$ & if $m\neq0$%
\end{tabular}
\\
&  =4^{-( n-1) }\sum_{j=0}^{k+1}2^{-j}+\biggl\{
\begin{tabular}
[c]{ll}%
$\dfrac{1}{\varepsilon}-3\cdot4^{-( n+k+1) }$ & if $m=0$\\
$4^{-m}$ & if $m\neq0$.
\end{tabular}
\end{align*}
\end{proof}

The case $m=0$ gives the next corollary.

\begin{corollary}
\label{C7}The set $Z(p_{n},p_{n+1}...,p_{n+k})$ has norm bounded by
$2\cdot4^{-(n-1)}+ {1/\varepsilon}.$
\end{corollary}

\begin{proposition}
\label{P8}Let $x$ be a vector in $Z(p_{n},...,p_{n+k})$, where $n\in
\mathbb{N}$ and $k\in\mathbb{N\cup}\{0\}.$ Then there exists a sequence of
pairwise disjoint vectors $(y_{j})_{j=0}^{k}$ such that
\[
x=\sum_{j=0}^{k}y_{j},\text{ }\Vert y_{j}\Vert_{\ell^{1}}=\dfrac{1}%
{\theta_{p_{n+j}}},\text{ }0\leq j\leq k
\]
and
\[
\operatorname*{supp}y_{j}\in\lbrack\mathcal{S}_{\eta_{n}+1},...,\mathcal{S}%
_{\eta_{n+j}+1}]\cap\lbrack M_{n+j}]^{<\infty}.
\]
\end{proposition}

\begin{proof}
The proof is by induction on $k.$ If $k=0,$ set $y_{0}=x$ and the claim is
clear. Assume the proposition holds for some $k$ and consider a vector $x\in
Z(p_{n},...,p_{n+k+1}).$ Write $x=\theta_{p_{n}}^{-1}\sum_{i=1}^{r}%
a_{i}e_{m_{i}}+\sum_{i=1}^{r}a_{i}z_{i}$ according to the definition of
$Z(p_{n},...,p_{n+k+1}).$\ By the inductive hypothesis, for each $i,$ there is
a sequence of pairwise disjoint vectors $(y_{j}^{i})_{j=1}^{k+1}$ such that
$z_{i}=\sum_{j=1}^{k+1}y_{j}^{i},$ $||y_{j}^{i}||_{\ell^{1}}=\theta_{p_{n+j}%
}^{-1},$ $\operatorname*{supp}y_{j}^{i}\in\lbrack\mathcal{S}_{\eta_{n+1+}%
1},...,\mathcal{S}_{\eta_{n+j}+1}]\cap\lbrack M_{n+j}]^{<\infty},$ $1\leq
j\leq k+1.$ Set $y_{0}=\theta_{p_{n}}^{-1}\sum_{i=1}^{r}a_{i}e_{m_{i}},$ and
$y_{j}=\sum_{i=1}^{r}a_{i}y_{j}^{i},$ $1\leq j\leq k+1.$ Then $(y_{j}%
)_{j=0}^{k+1}$ is a pairwise disjoint sequence such that $\sum_{j=0}%
^{k+1}y_{k}=x.$ Clearly, $\Vert y_{j}\Vert_{\ell^{1}}=\sum_{i=1}^{r}%
|a_{i}|||y_{j}^{i}||_{\ell^{1}}=\theta_{p_{n+j}}^{-1},$ $\ 1\leq j\leq k+1,$
and $\Vert y_{0}\Vert_{\ell^{1}}=\theta_{p_{n}}^{-1}\sum_{i=1}^{r}%
|a_{i}|=\theta_{p_{n}}^{-1}.$ Also, $\operatorname*{supp}y_{0}\in
\mathcal{S}_{\eta_{n}+1}\cap\lbrack M_{n}]^{<\infty}$ since $y_{0}\in
Z(p_{n}).$ Furthermore, since $m_{1}<y_{j}^{1}<...<m_{r}<y_{j}^{r}$ and
$\{m_{1},...,m_{r}\}\in\mathcal{S}_{\eta_{n}+1},$ $\operatorname*{supp}%
y_{j}\in\lbrack\mathcal{S}_{\eta_{n}+1},...,\mathcal{S}_{\eta_{n+j}+1}%
]\cap\lbrack M_{n+j}]^{<\infty},$ $1\leq j\leq k+1.$
\end{proof}


\begin{proof}
[Proof of Theorem \ref{T2}]Beginning with $M_{0}=M,$ carry out the
construction above. Now take a block basis $( z_{k}) $ of $( e_{n}) _{n\in
M\text{ }}$ such that $z_{k}\in Z( p_{1},p_{2},...,p_{k}) $ for all $k.$ By
Corollary \ref{C7}, $\| z_{k}\| \leq2+1/\varepsilon$ for all $k.$ Suppose
$F\in\mathcal{S}_{\omega^{\xi}}.$ Then there exists $j_{0}\leq\min F$ such
that $F\in\mathcal{S}_{\beta_{j_{0}}}.$ By Proposition \ref{P8}, for all $k\in
F,$ there exists $y_{k}$ such that $| y_{k}| \leq| z_{k}| ,$ $\| y_{k}\|
_{\ell^{1}}={\theta^{-1}_{p_{j_{0}}}}$ and $\operatorname*{supp}y_{k}%
\in\lbrack\mathcal{S}_{\eta_{1}+1},...,\mathcal{S}_{\eta_{j_{0}}+1}]\cap[
M_{j_{0}}] ^{<\infty}.$ Thus for all scalars $( a_{k}) ,$%
\[
\| \sum_{k\in F}a_{k}z_{k}\| \geq\| \sum_{k\in F}a_{k}y_{k}\| \geq
\theta_{p_{j_{0}}}\| \sum_{k\in F}a_{k}y_{k}\| _{\mathcal{F}_{p_{j_{0}}}}
=\theta_{p_{j_{0}}}\| \sum_{k\in F}a_{k}y_{k}\| _{\mathcal{\ell}^{1}},
\]
as $\mathcal{S}_{\beta_{j_{0}}}[ \mathcal{S}_{\eta_{1}+1},...,\mathcal{S}%
_{\eta_{j_{0}}+1}] \cap[ M_{j_{0}}] ^{<\infty}\subseteq\mathcal{F}_{p_{j_{0}}%
}.$ Therefore,%
\[
\| \sum_{k\in F}a_{k}z_{k}\| \geq\theta_{p_{j_{0}}}\sum_{k\in F}| a_{k}| \|
y_{k}\| _{\mathcal{\ell}^{1}}=\sum_{k\in F}| a_{k}| .
\]
\end{proof}

In the rest of the section, we prove the converse to Theorem \ref{T2}. By
\cite[Proposition 1]{LT}, we may assume without loss of generality that there
exists a sequence $(\ell_{n})\subseteq\mathbb{N}$ converging to $\infty$ such
that $\mathcal{F}_{n}=(\mathcal{F}_{n}\cap\lbrack\mathbb{N}_{\ell_{n}%
}]^{<\infty})\cup\mathcal{S}_{0}$ for all $n\in\mathbb{N}$, where
$\mathbb{N}_{k}=\{n\in\mathbb{N}:n\geq k\}.$

\begin{lemma}
\label{L3}If $(\dagger)$ fails, then for all $\varepsilon>0$ and all
$M\in\lbrack\mathbb{N}],$ there exist $M^{\prime}\in\lbrack M]$ and a regular
family $\mathcal{H}$ containing $\mathcal{S}_{0}$, $\iota(\mathcal{H}%
)<\omega^{\omega^{\xi}}$, such that for all sufficiently large $m,$ there
exist $n_{1},...,n_{s}$ so that $\varepsilon\theta_{n_{1}}\cdots\theta_{n_{s}%
}>\theta_{m}$ and $\mathcal{F}_{m}\cap\lbrack M^{\prime}]^{<\infty}%
\subseteq\lbrack\mathcal{H},\mathcal{F}_{n_{1}},...,\mathcal{F}_{n_{s}}].$
\end{lemma}

\begin{proof}
Fix $\varepsilon>0.$ Since $(\dagger)$ fails, there exists $\beta<\omega^{\xi
}$ such that for all $m,$ $\gamma(\varepsilon,m)+2+\beta\geq\ell(\alpha_{m}).$
Therefore, for all large enough $m,$ say $m>m_{0},$ there exist $n_{1}%
,...,n_{s}$ such that $\varepsilon\theta_{n_{1}}\cdots\theta_{n_{s}}%
>\theta_{m}$ and $\ell(\alpha_{n_{s}}...\alpha_{n_{1}})+2+\beta\geq\ell
(\alpha_{m}).$ Let $\beta^{\prime}=2+\beta+1<\omega^{\xi}.$ Then $\ell
(\alpha_{m})<\ell(\alpha_{n_{s}}...\alpha_{n_{1}})+\beta^{^{\prime}}.$ Thus,
\[
\iota(\mathcal{F}_{m})<\iota(\mathcal{S}_{\beta^{\prime}}[\mathcal{F}_{n_{1}%
},...,\mathcal{F}_{n_{s}}]).
\]
By the remark after Theorem \ref{G}, for all $N\in\lbrack\mathbb{N}],$ there
exists $N^{\prime}\in\lbrack N]$ such that
\[
\mathcal{F}_{m}\cap\lbrack N^{\prime}]^{<\infty}\subseteq\mathcal{S}%
_{\beta^{\prime}}[\mathcal{F}_{n_{1}},...,\mathcal{F}_{n_{s}}].
\]
Given $M\in\lbrack\mathbb{N}],$ applying the above argument repeatedly, we
obtain infinite sets
\[
M\supseteq M_{1}\supseteq M_{2}\supseteq...\supseteq M_{k}\supseteq\dots
\]
such that for all $k\in\mathbb{N},$ there exist $n_{1},...,n_{s}$ (depending
on $k$) such that $\varepsilon\theta_{n_{1}}\cdots\theta_{n_{s}}>\theta
_{m_{0}+k}$ and $\mathcal{F}_{m_{0}+k}\cap\lbrack M_{k}]^{<\infty}%
\subseteq\mathcal{S}_{\beta^{\prime}}[\mathcal{F}_{n_{1}},...,\mathcal{F}%
_{n_{s}}].$ Choose $(m_{k})$ so that $m_{0}<m_{1}<m_{2}<...$ and $m_{k}\in
M_{k}$ for all $k\in\mathbb{N}.$ Let $M^{\prime}=(m_{k})_{k=1}^{\infty}.$ For
all $k\in\mathbb{N},$ define $\mathcal{B}_{k}=\{G:\ell_{m_{0}+k}\leq G,\text{
}|G|\leq m_{k}\}$ and $\mathcal{B}=\cup_{k=1}^{\infty}\mathcal{B}_{k}%
\cup\mathcal{S}_{0}.$ Let $\mathcal{H=}(\mathcal{B},\mathcal{S}_{\beta
^{\prime}}).$ Then $\mathcal{H}$ contains $\mathcal{S}_{0}$ and
\[
\iota(\mathcal{H})=\iota(\mathcal{B},\mathcal{S}_{\beta^{\prime}}%
)=\iota(\mathcal{S}_{\beta^{\prime}})+\iota(\mathcal{B})=\omega^{\beta
^{\prime}}+\omega<\omega^{\omega^{\xi}}.
\]
Consider a set $F\in\mathcal{F}_{m_{0}+k}\cap\lbrack M^{\prime}]^{<\infty}$
for some $k\in\mathbb{N}$. Write $F=F_{1}\cup F_{2},$ where $F_{1}%
=F\cap\lbrack1,m_{k})$ and $F_{2}=F\cap\lbrack m_{k},\infty).$ Since $F_{1}%
\in\mathcal{F}_{m_{0}+k}=(\mathcal{F}_{m_{0}+k}\cap\lbrack\mathbb{N}%
_{\ell_{m_{0}+k}}]^{<\infty})\cup\mathcal{S}_{0},$ either $F_{1}\in
\mathcal{S}_{0}\subseteq\mathcal{B}$ or $F_{1}\in\mathcal{F}_{m_{0}+k}%
\cap\lbrack\mathbb{N}_{\ell_{m_{0}+k}}]^{<\infty}$. In the latter case,
$\ell_{m_{0}+k}\leq F_{1}$ and $|F_{1}|\leq m_{k}$ and hence $F_{1}%
\in\mathcal{B}_{k}\subseteq\mathcal{B}.$ Also, $F_{2}\in\mathcal{F}_{m_{0}%
+k}\cap\lbrack M_{k}]^{<\infty}$ implies that there exist $n_{1},\dots,n_{s}$
such that $\varepsilon\theta_{n_{1}}\cdots\theta_{n_{s}}>\theta_{m_{0}+k}$ and
$F_{2}\in\mathcal{S}_{\beta^{\prime}}[\mathcal{F}_{n_{1}},...,\mathcal{F}%
_{n_{s}}]$. Therefore, $F\in(\mathcal{B},\mathcal{S}_{\beta^{\prime}%
})[\mathcal{F}_{n_{1}},...,\mathcal{F}_{n_{s}}]=\mathcal{H}[\mathcal{F}%
_{n_{1}},...,\mathcal{F}_{n_{s}}].$
\end{proof}

\begin{proposition}
\label{LTP14}\emph{\cite[Proposition 14]{LT}} Suppose for all $\varepsilon>0,$
there exist a regular family $\mathcal{G}_{\varepsilon}$ and $m_{0}%
\in\mathbb{N}$ such that for all $m\geq m_{0},$ there exist $n_{1},\dots
,n_{s}\in\mathbb{N}$ satisfying $\theta_{m}<\varepsilon\theta_{n_{1}}%
\dots\theta_{n_{s}}$ and $\mathcal{F}_{m}\subseteq[ \mathcal{G}_{\varepsilon
},\mathcal{F}_{n_{1}},\dots,\mathcal{F}_{n_{s}}] .$ Then%
\[
I_{b}( X) \leq\sup_{\varepsilon>0}\sup\limits_{n\in\mathbb{N}}[ \iota(
\mathcal{G}_{\varepsilon}) \cdot\alpha_{n}^{\omega}] .
\]
\end{proposition}

\begin{theorem}
\label{conv} Suppose that $( \dagger) $ fails, then for all $M\in[ \mathbb{N}]
,$ there exists $N\in[ M] $ such that%
\[
I( [ ( e_{k}) _{k\in N}] ) =\omega^{\omega^{\xi}}.
\]
In particular, $[ ( e_{k}) _{k\in N}] $ does not contain any $\ell^{1}%
$-$\mathcal{S}_{\omega^{\xi}}$-spreading model.
\end{theorem}

\begin{proof}
By Lemma \ref{L3}, there exist infinite sets $M\supseteq M_{1}\supseteq
...\supseteq M_{k}\supseteq...$ such that for all $i\in\mathbb{N},$ there
exists a regular family $\mathcal{H}_{i}$ containing $\mathcal{S}_{0}$,
$\iota(\mathcal{H}_{i})<\omega^{\omega^{\xi}}$, such that for all sufficiently
large $n,$ say $n\geq m_{0}(i),$ there exist $n_{1},...,n_{s}$ so that
$\theta_{n}<\theta_{n_{1}}\cdots\theta_{n_{s}}/i$ and $\mathcal{F}_{n}%
\cap\lbrack M_{i}]^{<\infty}\subseteq\lbrack\mathcal{H}_{i},\mathcal{F}%
_{n_{1}},...,\mathcal{F}_{n_{s}}].$ Choose $m_{1}<m_{2}<m_{3}<...$ such that
$m_{k}\in M_{k}$ and let $N=(m_{k}).$ Set $Y=[(e_{k})_{k\in N}].$ Note that
$Y=T[(\theta_{n},\mathcal{G}_{n})_{n=1}^{\infty}],$ where $G\in\mathcal{G}%
_{n}$ if and only if $\{m_{k}:k\in G\}\in\mathcal{F}_{n}.$

Suppose $\varepsilon>0$ is given. Pick $i\in\mathbb{N}$ such that ${1}%
/{i}<\varepsilon.$ Assume that $n\geq m_{0}(i)$ and $\ell_{n}\geq m_{i}.$ If
$G\in\mathcal{G}_{n},$ then $F=\{m_{k}:k\in G\}\in\mathcal{F}_{n}\cap\lbrack
N]^{<\infty}.$ Since $\mathcal{F}_{n}=(\mathcal{F}_{n}\cap\lbrack
\mathbb{N}_{\ell_{n}}]^{<\infty})\cup\mathcal{S}_{0}$, either $F\in
\mathcal{S}_{0}$ or $F\in\mathcal{F}_{n}\cap\lbrack\mathbb{N}_{\ell_{n}%
}]^{<\infty}$. In the latter case, $F\geq\ell_{n}\geq m_{i}$ and thus
$F\in\mathcal{F}_{n}\cap\lbrack M_{i}]^{<\infty}.$ Hence in either case,
$F\in\lbrack\mathcal{H}_{i},\mathcal{F}_{n_{1}},...,\mathcal{F}_{n_{s}}]$ for
some $n_{1},...,n_{s}$ such that $\theta_{n}<\varepsilon\theta_{n_{1}}%
\cdots\theta_{n_{s}}.$ Therefore,%
\[
\mathcal{G}_{n}\subseteq\lbrack\mathcal{J}_{i},\mathcal{G}_{n_{1}%
},...,\mathcal{G}_{n_{s}}],
\]
where $G\in\mathcal{J}_{i}$ if and only if $\{m_{k}:k\in G\}\in\mathcal{H}%
_{i}.$ Note that $\iota(\mathcal{J}_{i})<\omega^{\omega^{\xi}}.$ Thus,
according to Proposition \ref{LTP14},%
\[
I_{b}(Y)\leq\sup_{i}\sup_{n\in\mathbb{N}}[\iota(\mathcal{J}_{i})\cdot
\alpha_{n}^{\omega}]=\sup_{i}[\iota(\mathcal{J}_{i})\cdot\omega^{\omega^{\xi}%
}]=\omega^{\omega^{\xi}}.
\]
However, $I_{b}(Y)\geq\omega^{\omega^{\xi}}$ by part 1 of \cite[Theorem
14]{LT}. Hence $I_{b}(Y)=\omega^{\omega^{\xi}}$. Finally, $I_{b}(Y)=I(Y)$ by
\cite[Corollary 5.13]{JO} since $I_{b}(Y)\geq\omega^{\omega}$. By \cite[Lemma
5.11]{JO}, $I(Y,K)<\omega^{\omega^{\xi}}$. Thus $Y$ does not contain any
$\ell^{1}$-$\mathcal{S}_{\omega^{\xi}}$-spreading model.
\end{proof}

\section{Mixed Tsirelson spaces constructed with Schreier families}

In this section, we apply the results of the last section to mixed Tsirelson
spaces of the type $T[ ( \theta_{n},\mathcal{S}_{\beta_{n}}) _{n=1}^{\infty}]
,$ where $( \theta_{n}) $ is a nonincreasing null sequence in $( 0,1) ,$
$\sup_{n}\beta_{n}=\omega^{\xi}>\beta_{n}>0$ for all $n\in\mathbb{N},$ and
$0<\xi<\omega_{1}.$ In the present situation, the function $\gamma$ is given
by
\[
\gamma(\varepsilon,m) = \max\{\beta_{n_{s}}+\cdots+\beta_{n_{1}}:
\varepsilon\theta_{n_{s}}\cdots\theta_{n_{1}} > \theta_{m}\}\ (\max\emptyset=
0).
\]
Theorems \ref{T2} and \ref{conv} give

\begin{theorem}
\label{sch} Let $(\beta_{n})$ be as above and let $(e_{n})$ be the unit vector
basis of the mixed Tsirelson space $T[(\theta_{n},\mathcal{S}_{\beta_{n}%
})_{n=1}^{\infty}]$. If condition \emph{($\dagger$)} holds, then for any
$M\in\lbrack\mathbb{N}]$, $(e_{n})_{n\in M}$ contains an $\ell^{1}%
$-$\mathcal{S}_{\omega^{\xi}}$-spreading model. If condition \emph{($\dagger
$)} fails, then for all $M\in\lbrack\mathbb{N}]$, there exists $N\in\lbrack
M]$ such that $[(e_{k})_{k\in N}]$ does not contain any $\ell^{1}%
$-$\mathcal{S}_{\omega^{\xi}}$-spreading model.
\end{theorem}

In the event that the Schreier families $\mathcal{S}_{\beta},$ $\beta$ a limit
ordinal, are defined using special choices, the second part of Theorem
\ref{sch} can be strengthened. The special ``standard'' choices are described
as follows. For all limit ordinals $\alpha<\omega_{1},$ fix a sequence of
ordinals strictly increasing to $\alpha.$ If $\beta=\omega^{\beta_{1}}\cdot
m_{1}+\dots+\omega^{\beta_{k}}\cdot m_{k}$ is a limit ordinal, determine
$\mathcal{S}_{\beta}$ using the sequence%
\[
\hat{\beta}_{n}=\bigl\{%
\begin{array}
[c]{ll}%
\omega^{\beta_{1}}\cdot m_{1}+\dots+\omega^{\beta_{k}}\cdot(m_{k}%
-1)+\omega^{\beta_{k}-1}\cdot n & \text{if }\beta_{k}\text{ is a successor}\\
\omega^{\beta_{1}}\cdot m_{1}+\dots+\omega^{\beta_{k}}\cdot(m_{k}%
-1)+\omega^{\zeta_{n}} & \text{if }\beta_{k}\text{ is a limit.}%
\end{array}
\]
where $(\zeta_{n})$ is the chosen sequence of ordinals increasing to
$\beta_{k}.$

\begin{theorem}
\emph{\cite[Theorem 26]{LT}} \label{LTT27A}Follow the notation above and apply
the standard choices to define Schreier families. If there exists
$\varepsilon>0$ such that for all $\beta<\omega^{\xi},$ there exists
$m\in\mathbb{N}$ satisfying $\gamma( \varepsilon,m) +2+\beta<\beta_{m},$ then
$I_{b}( T[ \theta_{n},\mathcal{S}_{\beta_{n}}) _{n=1}^{\infty}] )
=\omega^{\omega^{\xi}\cdot2}$. Otherwise, $I_{b}( T( \mathcal{F}_{0},(
\theta_{n},\mathcal{S}_{\beta_{n}}) _{n=1}^{\infty}) ) =\omega^{\omega^{\xi}}$.
\end{theorem}

For \textquotedblleft standard\textquotedblright\ Schreier families, the
second part of Theorem \ref{sch} can be improved.

\begin{theorem}
Let $(\beta_{n})$ be as above and apply the standard choices to define
Schreier families. If \emph{($\dagger$)} fails, then $I(T[(\theta
_{n},\mathcal{S}_{\beta_{n}}) _{n=1}^{\infty}]) = \omega^{\omega^{\xi}}$. In
particular, $T[(\theta_{n},\mathcal{S}_{\beta_{n}})_{n=1}^{\infty}]$ does not
contain any $\ell^{1}$-$\mathcal{S}_{\omega^{\xi}}$-spreading model.
\end{theorem}

Note that for finite $\beta_{n}$'s, no choices need to be made in defining the
Schreier families $\mathcal{S}_{n}.$ It is worthwhile to record the result in
this case.

\begin{theorem}
If $\theta_{m+n} \geq\theta_{m}\theta_{n}$ for all $m, n$ and $\lim_{m}%
\limsup_{n}\theta_{m+n}/\theta_{n}>0$, then $[ ( e_{k_{n}}) ] $ contains an
$\ell^{1}$-$\mathcal{S}_{\omega}$-spreading model for any subsequence $(
e_{k_{n}}) $ of the unit vector basis $( e_{k}) $ of $T[(\theta_{n}%
,\mathcal{S}_{n})_{n=1}^{\infty}]$. Otherwise $T[(\theta_{n},\mathcal{S}%
_{n})_{n=1}^{\infty}]$ contains no $\ell^{1}$-$\mathcal{S}_{\omega}$-spreading model.
\end{theorem}

\noindent\emph{Remark.} It can be shown that for sequences $(\theta_{n})$ such
that $\theta_{m+n} \geq\theta_{m}\theta_{n}$ for all $m, n$, the condition
$\lim_{m}\limsup_{n}\theta_{m+n}/\theta_{n}>0$ is strictly weaker than the
condition $\lim\theta^{1/n}_{n} = 1$.

\end{document}